\documentclass{amsart}
\usepackage{graphicx,amssymb,amscd,dsfont,psfrag}

\newtheorem{thm}{Theorem}[section]
\newtheorem{cor}[thm]{Corollary}
\newtheorem{lem}[thm]{Lemma}
\newtheorem{prop}[thm]{Proposition}
\newtheorem*{mthm}{Main Theorem}
\theoremstyle{definition}
\newtheorem{defn}[thm]{Definition}
\theoremstyle{remark}
\newtheorem{rem}[thm]{Remark}

\newtheorem*{question}{Question}
\newtheorem*{example}{Example}
\numberwithin{equation}{section}

\newcommand{\set}[1]{\left\{#1\right\}}

\newcommand{\To}{\longrightarrow}

\newcommand{\A}{\mathcal{A}}
\newcommand{\D}{\mathcal{D}}
\newcommand{\W}{\mathcal{W}}
\newcommand{\Y}{\mathsf{Y}}
\newcommand{\p}{\underline{p}}

\newcommand{\Apc}{\A_{c}(p)_{-1}}
\newcommand{\Dpc}{\D_{c}(p)_{-1}}
\newcommand{\nvct}{(n_{1},n_{2},\ldots,n_{p})}

\DeclareMathOperator{\Wltemp}{\W}
\newcommand{\Wl}{\sideset{^{l}}{}{\Wltemp}}
\DeclareMathOperator{\Hltemp}{H}
\newcommand{\Hl}{\sideset{^{l}}{}{\Hltemp}}

\DeclareMathOperator{\Altemp}{\A}

\newcommand{\Apcl}{\sideset{^{l}}{_{c}}{\Altemp}(p)_{-1}}
\DeclareMathOperator{\Dltemp}{\D}

\newcommand{\Dpcl}{\sideset{^{l}}{_{c}}{\Dltemp}(p)_{-1}}
\def\co{\colon\thinspace}
\def\ass{{\text{\rm \raisebox{0.03ex}{:}\vspace{-0.05ex}=}}}

\begin{document}

\title[Acyclic Jacobi Diagrams]
{Acyclic Jacobi Diagrams}%
\author{Daniel Moskovich}%
\address{Research Institute for Mathematical Sciences,
Kyoto University, Kyoto, 606-8502 JAPAN}%
\email{dmoskovich@gmail.com}%
\urladdr{http://www.sumamathematica.com/}

\thanks{I would like to thank Tomotada Ohtsuki for his perennial help and support,
Kazuo Habiro for useful comments, Fred Cohen for some stimulating
discussions, and Atsushi Ishii, Tadayuki Watanabe, and Eri
Hatakenaka for being there for me when I was in hospital.
Special thanks also to Sergei Duzhin for his careful reading and corrections, and to the anonymous referee for many useful
comments.}%
\subjclass{57M27,17B01}%
\keywords{Free Lie algebra, Jacobi diagram, Vassiliev invariants,
Milnor invariants, link invariants, left-normed basis}%

\date{6th of September, 2006}%
\begin{abstract}
We propose a simple new combinatorial model to study spaces of
acyclic Jacobi diagrams, in which they are identified with algebras
of words modulo operations. This provides a starting point for a
word-problem type combinatorial investigation of such spaces, and
provides fresh insights on known results.
\end{abstract}
\maketitle
\section{Introduction}

Jacobi diagrams are a subset of labeled pseudo-graphs whose vertices
have valence $1$ or $3$, with some extra structure. They provide a
profound and as yet largely unexplored bridge between the world of
Lie algebras and the world of low-dimensional topology, concisely
encoding topological invariants which are in some sense ``Lie
algebra-like''.\par

We work over a field $\Bbbk$ of characteristic different from $2$.
Let $\mathds{N}_{0}$ denote the non-negative integers (the natural
numbers). Roughly speaking, to a fixed Lie algebra (or more
generally to a fixed Lie algebra object \cite{NW04}), and a
low-dimensional topological object $M$ (a knot, a link, a
$3$--manifold, a hyper-K\"{a}hler manifold, a handlebody\ldots), a
Jacobi diagram defines a map called a \emph{weight system} from the
input data to an $\mathds{N}_{0}$--graded vector space $V$. Thus,
Jacobi diagrams plus weight systems define $V$--valued topological
invariants called the \emph{finite type invariants} of $M$.\par

If the object of study is the space of finite type invariants
(either of a given $M$ or of all $M$), the key space to understand
is the space of Jacobi diagrams (see \cite{NW04} for the
set-theoretical arguments allowing us to call this a space or a
set). In particular, being able to identify this space with a space
of Lie algebras (or of Lie algebra objects) with something cut away
would be a huge step forward both for low dimensional topology and
for the theory of Lie algebras.\par

For the space of all Jacobi diagrams the question above is wide open
(see however \cite{CKX03} and \cite{Koh97} for recent progress). But
when we restrict to acyclic graphs--- to Jacobi diagrams without
loops--- we have a some understanding of the relationship between
the Lie algebra world and the world of Jacobi diagrams
(\cite[Section 4.3]{Mor93}, \cite{CohF05,Lev02,Lev06}). We should
note however that neither side of this relationship is in itself
well understood.

\begin{thm}\label{T:hn}
The space of connected acyclic Jacobi diagrams over $\mathds{Q}$ is
isomorphic to $\mathbf{h}_{n}(\p)$ which is given by
\begin{equation}
\mathbf{h}_{n}(\p)\To \mathcal{L}_{n-1}(\p)\otimes\p\To
\mathcal{L}_{n}(\p) \To 0
\end{equation}
When $\mathcal{L}_{n}(\p)$ denotes the length $n$ part of the free
Lie algebra over $p$ generators.\footnote{See \cite{Lev02,Lev06} for
the corresponding question over $\mathds{Z}$.}
\end{thm}

One corollary of our main theorem is a simple combinatorial proof of
the above result, first proved in the $\Bbbk=\mathds{Q}$ case in
\cite{HP03,Lev02}. 

An example of Theorem \ref{T:hn} for $n=p=3$ is given below in
set--theoretical notation

\begin{multline*}
\left\{\hspace{7pt}\begin{minipage}{80pt}
\psfrag{a}{$1$}\psfrag{b}{$2$}\psfrag{c}{$3$}
\includegraphics[width=80pt]{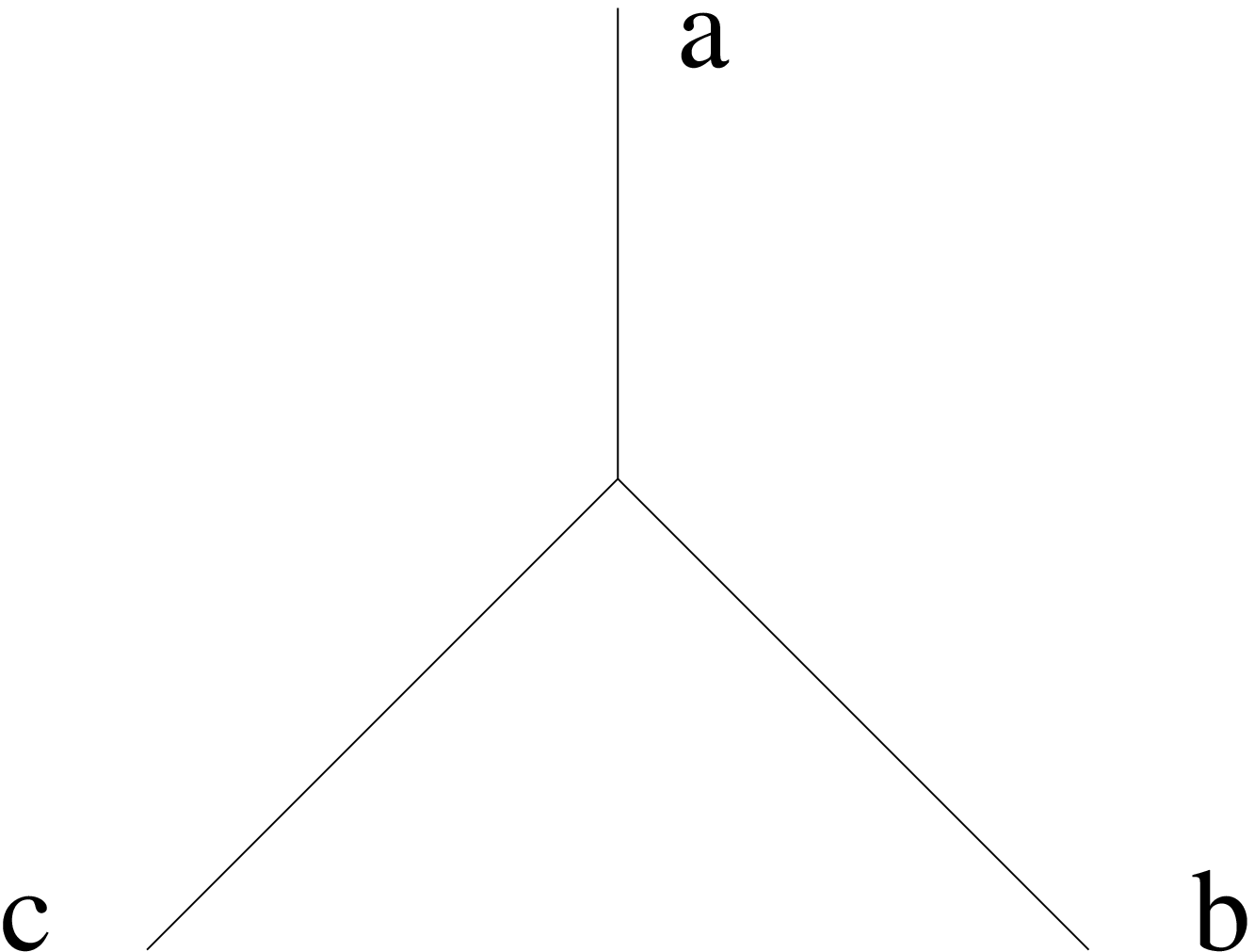}
\end{minipage}\hspace{7pt}\right\}\To
\begin{Bmatrix}
  \ [1,2]\otimes 1 & [1,2]\otimes 2 & [1,2]\otimes 3\ \\
  \ [1,3]\otimes 1 & [1,3]\otimes 2 & [1,3]\otimes 3\  \\
  \ [2,3]\otimes 1 & [2,3]\otimes 2 & [2,3]\otimes 3\
\end{Bmatrix}\\ \To
\begin{Bmatrix}
\ [[1,2],1] & [[1,2],2] & [[1,2],3]\  \\
\ [[1,3],1] & \text{---} & [[1,3],3]\  \\
\ [[2,3],1] & [[2,3],2] & [[2,3],3]\
\end{Bmatrix}
\To 0
\end{multline*}


The topological meaning of the finite type invariants associated to
acyclic Jacobi diagrams is well known--- these are the Milnor
$\bar{\mu}$--invariants which measure linkage
(\cite{HL98,HM00,Pol99}). These invariants play no front-line role
in the study of knots and of manifolds (although they play an
important role behind the scene as part of what is known as the
Associator), but they come to the fore in the study of finite-type
invariants of links, braids, tangles.\par

In this paper we show that the space of acyclic Jacobi diagrams is
equivalent to a subspace of itself of graphs of a certain special
shape called \emph{swings} modulo a set of moves inherited from the
moves on Jacobi diagrams. The space of swings is computationally
simpler than the full space of acyclic Jacobi diagrams, and is
isomorphic either to an algebra of words we call $\W'(p)$ or to a
different algebra of words we call $\Wl(p)$. These are defined to be
the free associative $\Bbbk$--algebra on $p$ letters $ASS(\p)$
modulo the families of fold moves $H'$ and $\Hl$ correspondingly
(see Section \ref{SS:sw}). Which of these two algebras of words the
space of swings is isomorphic to depends on whether we are dealing
with unrooted or with rooted Jacobi diagrams (free Lie
algebras).\par

Formally stated:

\begin{mthm}
The space of (unrooted) acyclic Jacobi diagrams $\Apc$ is isomorphic
to the algebra of words $\W'(p)$ and the free Lie algebra
$\mathcal{L}(\p)$ is isomorphic to the algebra $\Wl(p)$.
\end{mthm}

An example of the first part of the theorem is given below

$$
\left.\left\{\left.\hspace{7pt}
\psfrag{a}{$a$}\psfrag{b}{$b$}\psfrag{c}{$c$}\psfrag{d}{$d$}\psfrag{e}{$e$}\psfrag{f}{$f$}
\begin{minipage}{40pt}
\includegraphics[width=40pt]{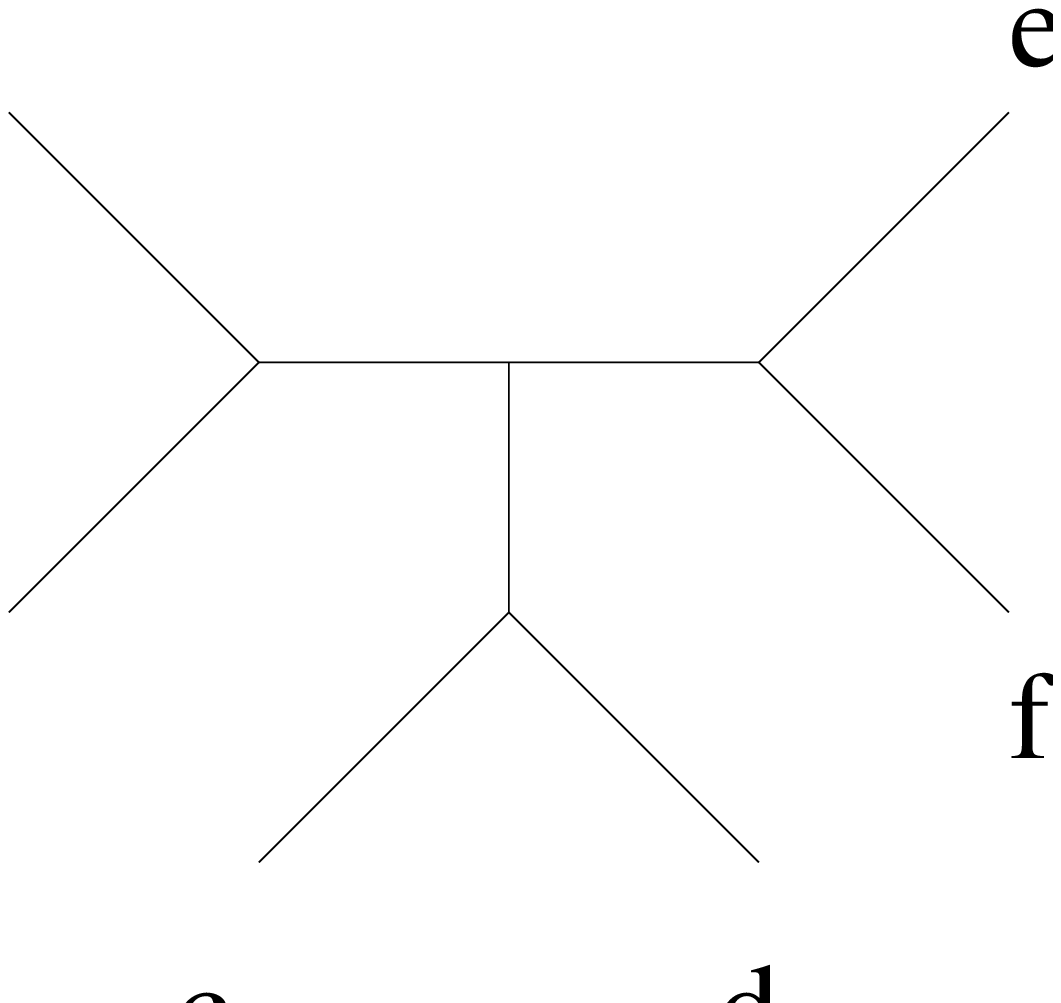}
\end{minipage}\hspace{12pt}\raisebox{-10pt}{,}\hspace{5pt}
\begin{minipage}{60pt}
\includegraphics[width=60pt]{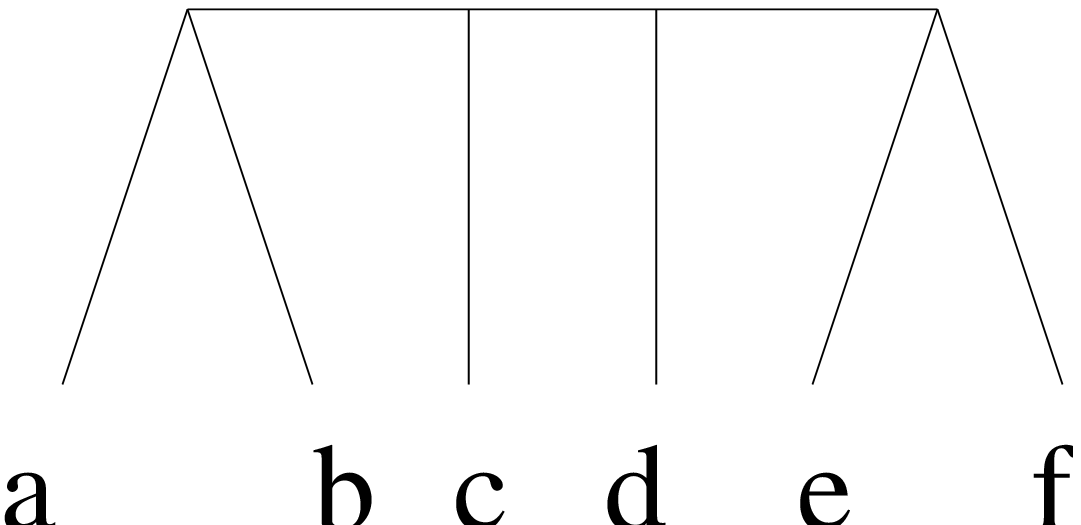}
\end{minipage}\hspace{10pt}\right|\ a,b,c,d,e,f\in\underline{6}
\right\}\right/AS\text{, }IHX\simeq ASS(\underline{6})/H'
$$

Definitions of the above concepts will be given in Section
\ref{S:defn}. The space of connected acyclic Jacobi diagram chains
$\Apc$ and its rooted version $\mathcal{L}(\p)$ the free Lie algebra
over $\Bbbk$ will be defined in Section \ref{SS:ajd}. The remaining
terminology and notation which pertains to spaces of swings and
their moves will be defined in Section \ref{SS:sw}.\par

In Section \ref{S:sample} we illustrate how our main theorem may be
used to calculate by hand the $5373540$ basis elements of
$\mathbf{h}_{9}(\underline{9})$. This example provides some basis
for speculation on the form of a swing-basis for general spaces of
acyclic Jacobi diagrams.\par

We conclude this introduction with a short summary of the current
state of knowledge about $\mathbf{h}_{n}(\p)$ and how our result
supplements it. Theorem \ref{T:hn} gives dimensions of spaces of
acyclic Jacobi diagrams in terms of dimensions of free Lie algebras.
These are given by Witt's dimension formula \cite{Wit37}:

$$
\dim(\mathcal{L}_{n}(\p))=\frac{1}{n}\sum_{d\mid
n}\mu(d)p^{\frac{n}{d}}
$$

The sum is over all (positive) divisors $d$ of $n$. The M\"{o}bius
function $\mu(d)$ is defined by

$$
\mu(d)\co\left\{
  \begin{array}{ll}
    1, & \hbox{If $d=1$;} \\
    (-1)^{k}, & \hbox{If $d=p_{1}\cdots p_{k}$ (distinct $p_{i}$);} \\
    0, & \hbox{if $d$ has a square factor.}
  \end{array}
\right.
$$

Let now $\mathcal{L}(\p)_{\nvct}$ denote the subalgebra of the free
Lie algebra over $p$ letters consisting of words in which the $i$th
element appears $n_{i}$ times for $1\leq i\leq p$, where
$\sum_{i=1}^{p}n_{i}=n$. Then the dimension of this space, the
\emph{necklace number}, is given by the following formula also due
to Witt \cite{Wit37}:
$$
\dim(\mathcal{L}(\p)_{\nvct})=\frac{1}{n}\sum_{d\mid
n_{1},\ldots,n_{p}}\mu(d)\frac{\frac{n}{d}!}{\frac{n_{1}}{d}!\cdots\frac{n_{p}}{d}!}
$$



To make further computational progress for acyclic Jacobi diagrams,
the best basis one might hope for would be a monomial basis for
$\mathbf{h}_{n}(\p)$ in terms of swings. Theorem \ref{T:hn} reduces
the problem to that of finding a monomial basis for free Lie
algebras. Sergei Duzhin suggests that this might be an important
step in the calculation of the rational associator in general, not
just when restricting to acyclic Jacobi diagrams.\par

The problem of calculating left-normed bases for free Lie algebras
was first discussed by Kukin \cite{Kuk78} who claimed to have solved
it. Fifteen years later, deficiencies in his construction were
revealed by Blessenohl and Laue \cite{BlL93} who offered an
alternative construction which works over any field $\Bbbk$ which
contains all roots of unity. Unfortunately this condition rarely
holds in a topological setting, where the ground field is usually
the rationals, the integers, or some finite field. The main theorem
of this paper significantly simplifies the algorithmic calculation
of such a basis. Since it appears unlikely that our set of moves on
swings is minimal, this calculation can likely be simplified yet
further. We hope to return to this problem in the future.\par


This paper is a reorganized version on the sections on acyclic
Jacobi diagrams in the author's Master Thesis in the University of
Tokyo \cite{Mos03}.

\section{Basic Definitions and Notation}\label{S:defn}

\subsection{Acyclic Jacobi Diagrams}\label{SS:ajd}

Fix a natural number $p$ and a field $\Bbbk$ of characteristic
different from $2$.

\begin{defn}
An acyclic Jacobi diagram is a connected vertex-oriented acyclic
graph whose vertices have valence $1$ or $3$ and whose univalent
vertices (legs) are labeled by elements of $\p$ the ordered set
$\set{1,2,\ldots,p}$ (viewed as a vector space with ordered basis so
that we can tensor it with vector spaces).
\end{defn}

All graphs discussed in this paper come equipped with a fixed
arbitrary ordering of their vertices and arcs.\par

Concepts defined for acyclic graphs naturally specialize to the case
of acyclic Jacobi diagrams. In the standard terminology for the
Jacobi diagram world, a univalent vertex $s$ of a graph $G$ is
called a \emph{leg} of $G$ \cite{BN95b}.\par

Two classes of graphs with specified vertex will be used in this
paper. A \emph{pointed graph} $\sigma_{s}(G)$ is a graph $G$ with a
distinguished leg $s$ \cite{Joy81,BLRL97}. A \emph{rooted Jacobi
diagram} $rt_{s}(G)$ is a Jacobi diagram $G$ with a single leg $s$
labeled by a distinguished element $*$ called the \emph{root} rather
than by an element of $\p$ \cite{Ber85}. When $s$ is clear from the
context we may omit it from the notation and we may denote such
graphs $\sigma(G)$ and $rt(G)$ correspondingly.\par

Rooted acyclic Jacobi diagrams (acyclic Jacobi diagrams with a
single leg labeled $*$) are in bijective correspondence with $M(\p)$
the free magma over $p$ letters. The binary operation is `connecting
at the root' \cite{Ser66}.\par

Returning to Jacobi diagram specific terminology, an \emph{acyclic
Jacobi diagram chain} is a linear sum of acyclic Jacobi diagrams
over $\Bbbk$. This space is denoted $\Dpc$. Pointing $\Dpc$ gives us
the algebra $\Bbbk(M(\p))$.

\begin{rem}
The notation $\Dpc$ is an attempt to be consistent with the standard
notation for Jacobi diagrams in general, as outlined in
\cite{BN95b}. In the present paper we shall not be interested in
what the subscripts, superscripts, and $-1$ in the bracket denote.
$\Dpc$ is to be taken as though it were a single compound symbol.
\end{rem}

The space of acyclic Jacobi diagrams comes equipped with two
operations, which are local moves between acyclic Jacobi diagrams
which have embeddings into $\mathds{R}^{2}$ which differ inside a
dotted circle as indicated below.
\begin{enumerate}
\item[AS]$$
\begin{minipage}{42pt}
    \includegraphics[width=42pt]{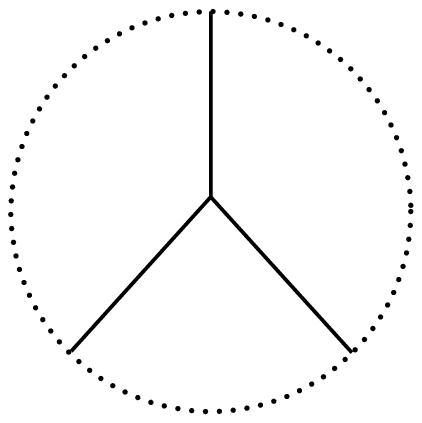}
\end{minipage}
\hspace{10pt}=\hspace{10pt}-\hspace{5pt}
\begin{minipage}{42pt}
    \includegraphics[width=42pt]{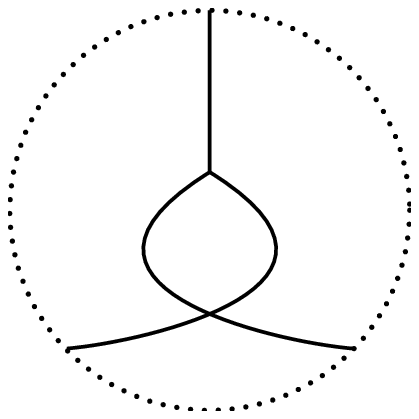}
\end{minipage}
$$

\item[IHX]$$
\begin{minipage}{42pt}
    \includegraphics[width=42pt]{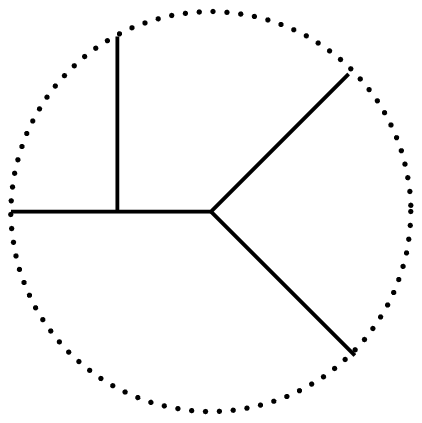}
\end{minipage}
\hspace{10pt}=\hspace{10pt}
\begin{minipage}{42pt}
    \includegraphics[width=42pt]{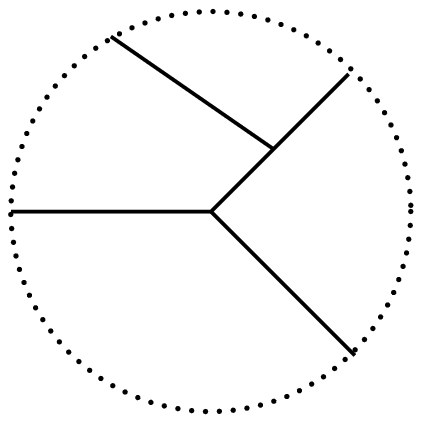}
\end{minipage}
\hspace{10pt}+\hspace{10pt}
\begin{minipage}{42pt}
    \includegraphics[width=42pt]{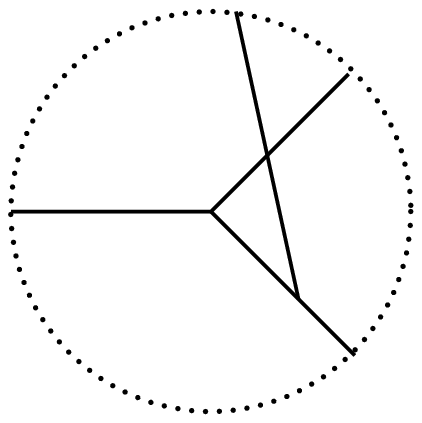}
\end{minipage}
$$
\end{enumerate}

Equivalence classes of acyclic Jacobi diagrams modulo $AS$ and $IHX$
are called \emph{acyclic Jacobi diagram classes}. This space is
denoted $\Apc$.

\begin{rem}
In the literature, Jacobi diagrams, Jacobi diagram chains, and
Jacobi diagram classes are all called Jacobi diagrams (except in
\cite{NW04} in which the distinction is made between the first two
concepts). In the present paper the distinction between these terms
becomes important because of our direct combinatorial approach. From
the point of view of notation however we will not distinguish
between acyclic Jacobi diagrams and acyclic Jacobi diagram classes.
\end{rem}

The rooted version of $\Apc$ is $\mathcal{L}(\p)$, the free Lie
algebra on $p$ generators. $AS$ corresponds to the anti-symmetry
relation, which $IHX$ is the Jacobi identity (this is the reason for
the name `Jacobi diagram' \cite{BN95c}).\par

The spaces above have a natural degree grading given by the number
of legs in the graphs. Denote this degree by $len(G)$, the
\emph{length} of $G$. When a finer grading is needed, we may refine
the length to get $mld(G)$ the \emph{multi-degree} of a graph $G$.
This is defined to be the vector $\nvct$ where $n_{i}$ is the number
of legs of $G$ coloured by the $i$th element of $\p$.

\subsection{Swings}\label{SS:sw}

In this section we introduce the space of words with which we would
like to model the space of acyclic Jacobi classes $\Apc$. This space
is called the space of \emph{swings}. As a set it is a subset of
$\Apc$, but its set of operations is smaller and more
computer-friendly.\par

A \emph{vertebrate} $G_{s_{1},s_{2}}$ is an acyclic uni-trivalent
graph $G$ with distinguished legs $s_{1}$ and $s_{2}$ called the
\emph{head} and the \emph{tail} of $G$ correspondingly. The
\emph{vertebral column} of $G$ is the unique elementary path from
the tail of $G$ to its head. If $s_{1}=s_{2}$ then $G_{s_{1},s_{2}}$
is said to be \emph{degenerate} \cite{Joy81,BLRL97}. There is a map
$\sigma_{s_{1}}\sigma_{s_{2}}$ from acyclic uni-trivalent graphs to
vertebrates that chooses $s_{1}$ as the head and $s_{2}$ as the tail
(this map is defined because the vertices are ordered. In Section
\ref{SS:headtail} we prove that in a suitable sense this map is
independent of this ordering.).\par

\begin{figure} \label{Fi:Breaktree}
\begin{center}
\psfrag{T1}{$T_{1}$}\psfrag{T2}{$T_{2}$}\psfrag{IHX}{$IHX$}
    \scalebox{.73}{\includegraphics{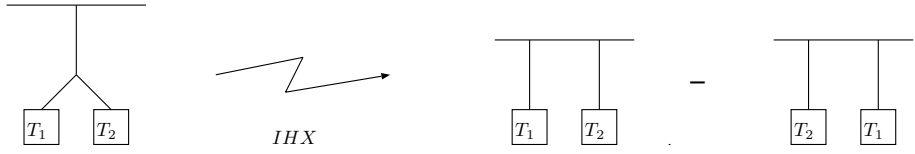}}
    \caption{Breaking a tree into branches}-
\end{center}
\end{figure}

A vertebrate $G_{s_{1},s_{2}}$ is said to be a \emph{swing} if all
trivalent vertices of $G_{s_{1},s_{2}}$ are neighbours of legs. The
`rooted version' $rt_{s_{1}}(G_{s_{1},s_{2}})$ is called a
\emph{half-swing}. Swings and half-swings are the main objects of
this paper, which we shall use to model the whole algebra of acyclic
Jacobi diagrams. Define $Sw(-)$ to be the restriction of a
$\Bbbk$--algebra of acyclic Jacobi diagrams to its
$\Bbbk$--subalgebra of $\Bbbk$--linear sums of swings, or of a
$\Bbbk$--algebra of rooted acyclic Jacobi diagrams to its
$\Bbbk$--subalgebra of $\Bbbk$--linear sums of half-swings.\par

In Section \ref{SS:rho} we construct a well--defined map $\rho\co
ASS(M(\p))\To Sw(\Dpc)$ from the free associative algebra generated
by elements of $M(\p)$. Correspondingly, there is a map $\rho^{l}$
from $\Apcl$ to $Sw(\mathcal{L}(\p))$.\par

A vertebrate may be thought of as an element of a free associative
algebra. Think of rooted trees (elements of $M(\p)$) as generators
and the degenerate rooted tree as the identity element, and read off
the product along the vertebral column from tail to head. In this
way, a swing corresponds to a word on $p$ letters. The free
associative $\Bbbk$--algebra of words on $p$ letters is denoted
$ASS(\p)$ and the free associative $\Bbbk$--algebra of words on a
set $\mathcal{X}$ of rooted trees is denoted $ASS(\mathcal{X})$. The
free associative product is given by concatenation of vertebrates,
gluing the head of one to the tail of the next.\par

We now define a number of actions on words in $ASS(\p)$. The actions
extend to chains by linearity.

\begin{defn}
\begin{equation}
\eta(\prod_{i=1}^{n}a_{i}):=
\begin{cases}
a_{1} &\text{if $n=1$;}\\
a_{n}\eta(\prod_{i=1}^{n-1}a_{i})-\eta(\prod_{i=1}^{n-1}a_{i})a_{n}
&\text{otherwise.}
\end{cases}
\end{equation}

$\eta$ is the algebraic operation corresponding to Figure 1 when
$T:=T_{1}T_{2}$ is a rooted tree corresponding to a left-bracketed
word in $M(\p)$.
\end{defn}

\begin{defn}
\hfill
\begin{itemize}

\item For a word $w\in ASS(\p)$:
\begin{equation}
h^{l}_{n}(w):=
    \begin{cases}
        (-1)^{n-1}a_{n}(\eta(\prod_{i=1}^{n-1}a_{i})(\prod_{i=n+1}^{len(w)}a_{i})
        &\text{for $2\leq n\leq len(w)$;}\\
        w   &\text{otherwise.}
\end{cases}
\end{equation}

\item For a word $w\in ASS(\p)$:
\begin{equation}
h'_{n}(w):=
\begin{cases}
0 &\text{for $len(w)=1$;}\\
 (-1)^{n}\prod_{i=1}^{n}a_{n-i+1} &\text{for
$n=len(w)>1$;} \\
h^{l}_{n} &\text{otherwise.}
\end{cases}
\end{equation}
\end{itemize}
\end{defn}

$h^{l}_{n}$ and $h'_{n}$ shall be called \emph{fold moves}, the
image being of folding egg-whites into a mixture. These are the
basic operations on swings in this paper.

\begin{figure} \label{Fi:haction}
\begin{center}
\psfrag{w1}{$w_{1}$}\psfrag{w2}{$w_{2}$}\psfrag{an}{$a_{n}$}\psfrag{IHX}{$IHX$}\psfrag{AS}{$AS$}
\psfrag{h}{$\eta$}\psfrag{wl}{$(w_{1})$}\psfrag{rotate}{rotate}
\scalebox{.77}{\includegraphics{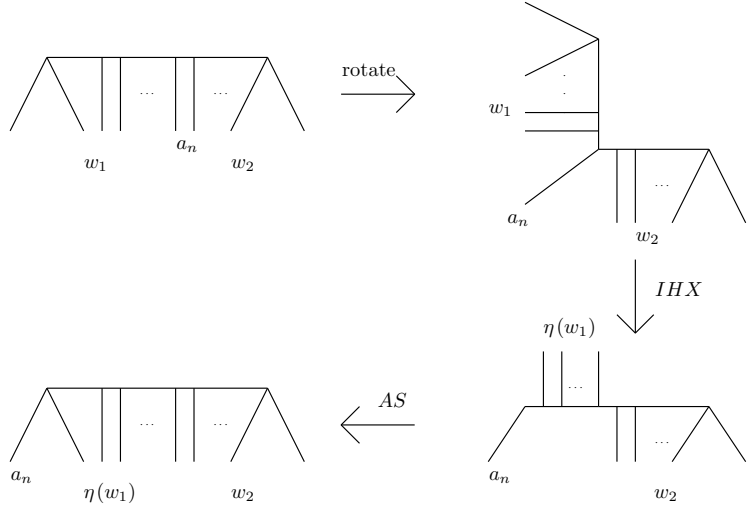}}
    \caption{The action on a swing that $h'_{n}$ is modeling. Here $w_{1}:=\prod_{i=1}^{n-1}a_{i}$ and
$w_{2}:=\prod_{i=n+1}^{len(w)}a_{i}$.}
\end{center}
\end{figure}

\begin{defn}
Let $\Hl$ be the set of all relations of the form $h^{l}_{n}(w)-w=0$
and let $H'$ be the set of all relations of the form $h'_{n}(w)-w=0$
for all $n\geq 2$ and for all $w\in ASS(\p)$. We refer to $\Hl$ and
to $H'$ collectively as the \textit{$H$--relations} or as the
collection of \emph{fold moves}. Define $\W'(p)$ to the quotient of
$ASS(\p)$ by $H'$, and define $\Wl(p)$ to the quotient of $ASS(\p)$
by $\Hl$.
\end{defn}

\begin{rem}
The naming $\Wl$ stems from the fact that $h^{l}_{n}$ maps act on
the left of the word. The superscript on the upper left hand side
serves to make sure that $\W'$ and $\Wl$ do not look too much alike.
\end{rem}

In the next section we shall use the actions we have defined here on
$ASS(\p)$ to recover well-known facts about free Lie algebras in our
language and setting. The basic reference for the material below is
Christophe Reutenauer's book \textsl{Free Lie Algebras}
\cite{Reu93}.

\subsection{Free Lie Algebra Identities in Terms of $\A$--Spaces.}

\subsubsection{Relations Involving $\eta$}

\begin{lem} [\cite{Reu93}, Theorem 1.4 (v)] \label{L:etaneta}
$$\eta(\eta(w))=(-1)^{n-1}n\eta(w)$$ for any $w\in ASS(\p)$.
\end{lem}

\begin{proof}
It is sufficient to prove the claim for words of multidegree
$(1,1,\ldots,1)$, and by linearity it is enough to prove it for $w$
a single word with coefficient $1$. Our proof shall be by induction.
For $len(w)=2$, we have
$\eta(\eta(w))=\eta(a_{2}a_{1}-a_{1}a_{2})=2\eta(w)$. Let us assume
that the claim is true until a certain length $n$. For $w$ of length
$n+1$, then, $n$ elements of $\eta(w)$ with sign $(-1)^{n-1}$ appear
in $\eta(\eta(w))$ by induction as $\eta$ of words ending in
$a_{n+1}$ in $\eta(w)$. They then appear again as reflections of
their reflections which are in $\eta(w)$, when the reflection
changes signs or preserves signs both times. These are all the
elements in $\eta(\eta(w))$ which start or end in $a_{n+1}$.\par

Now all we have to show is that there are no words in
$\eta(\eta(w))$ in which $a_{n+1}$ does not appear either in the
first place or in the last place. Let us assume that there is such a
word, of the form $w_{1}a_{n+1}w_{2}$. This word may come either
from $w_{1}\tau(w_{2})a_{n+1}$ or from $\tau(w_{2})w_{1}a_{n+1}$ (as
$a_{n+1}$ must be in the first or last place in $\eta(w)$). These
appear the same number of times in $\eta(w)$, with opposite signs---
when $n$ is odd, if $w_{1}\tau(w_{2})a_{n+1}$ appears with plus,
then $w_{2}\tau(w_{1})a_{n+1}$ would be plus, and then
$\tau(w_{2})w_{1}a_{n+1}$ comes out with a minus. When $n$ is even,
if $w_{1}\tau(w_{2})a_{n+1}$ appears with plus, then
$w_{2}\tau(w_{1})a_{n+1}$ gets a minus sign, and then
$\tau(w_{2})w_{1}a_{n+1}$ preserves the sign because $w_{1}$ and
$w_{2}$ are of the same parity. Then, from both
$w_{1}\tau(w_{2})a_{n+1}$ and $\tau(w_{2})w_{1}a_{n+1}$, the
sub-word $\tau(w_{2})$ is reflected, and every sub-word has the same
parity as itself, therefore the sign does not change, and the
elements cancel each other out.
\end{proof}

\begin{cor} \label{C:etaetakill}
$$\eta(\eta(w_{1})\eta(w_{2})+\eta(w_{2})\eta(w_{1}))=0$$
for any $w_{1},w_{2}\in ASS(\p)$ of lengths $n_{1},n_{2}$
correspondingly, $n_{1}+n_{2}>2$. In particular,
$\eta(\eta(w_{1})\eta(w_{1}))=0$.
\end{cor}

\begin{proof}
Precisely as in the proof of Lemma \ref{L:etaneta}, all elements of
$\eta(\eta(w_{1})\eta(w_{2}))$ in which a word in $\eta(w_{1})$ is
split into two parts are canceled, and we are left with
$(-1)^{n}\eta(\eta(w_{2})))\eta(\eta(w_{1}))+(-1)^{n-1}\eta(\eta(w_{1}))\eta(\eta(w_{2}))$
(second term: $(-1)^{n_{2}}$ from the $w_{2}$) where $n$ is defined
to be the combined length of $w_{1}$ and of $w_{2}$. By Lemma
\ref{L:etaneta}, this equals
$n(\eta(w_{2})\eta(w_{1})-\eta(w_{1})\eta(w_{2}))$. For
$\eta(\eta(w_{1})\eta(w_{2}))$, we similarly obtain
$n(\eta(w_{1})\eta(w_{2})-\eta(w_{2})\eta(w_{1}))$. The sum of these
terms is zero.
\end{proof}

\begin{lem}[\textit{An Identity of Baker (1905)}, \cite{Reu93} Section 1.6.6] \label{L:Y2I}
For any $w_{1},w_{2}$ words of lengths $n_{1}, n_{2}$ respectively
in $ASS(\p)$,
$$\eta(w_{1})\eta(w_{2})=(-1)^{n_{2}}(\eta(w_{1})\eta(w_{2})-\eta(w_{2})\eta(w_{1}))$$
\end{lem}

\begin{proof}
$(-1)^{n_{2}}(\eta(w_{1})\eta(w_{2})-\eta(w_{2})\eta(w_{1}))$ are
the elements of $\eta(w_{1}(\eta(w_{2})))$ where $\eta(w_{2})$ is
preserved. We must show that rest of the elements in the image of
$w_{1}(\eta(w_{2}))$ under $\eta$ cancel out. This follows from the
claim in the proof of Lemma \ref{L:etaneta}, that if an element of
$\eta(w_{2})$ is split into 2 parts, $w'$ and $w''$ so that we get
$w'w'_{1}w''$ where $w'_{1}$ is a word in $\eta(w_{1})$, such an
element can come only from $w'_{1}w''\tau(w')$ and from
$w'_{1}\tau(w')w''$, and such elements cancel out (same sign
argument as above precisely).
\end{proof}

\subsubsection{Relations Involving Fold Moves}

\begin{lem} \label{L:11}
For $w:=\prod_{i=1}^{n}a_{i}$, if $a_{i}=a_{i+1}$, then
$h_{i}^{l}(w)=h_{i+1}^{l}(w)$.
\end{lem}

\begin{proof}
Direct calculation.
\end{proof}

\begin{lem} \label{L:hreduce}
Let $w:=\prod_{i=1}^{n}a_{i}$ be a word in $ASS(\p)$. Then for
$i>j$, we have $h^{l}_{i}h^{l}_{j}(w)=h^{l}_{i}(w)$.
\end{lem}

\begin{proof}
Let $w'$ be $\prod_{i=1}^{j-1}a_{i}$. By Lemma \ref{L:etaneta}, we
have that
$$h^{l}_{i}(w)=\frac{(-1)^{j}}{j-1}h^{l}_{i}(\eta(w')\prod_{i=j}^{n}a_{i})$$
But now we have the equality
$$h^{l}_{i}(w)-\frac{(-1)^{j}}{j-1}h^{l}_{i}h^{l}_{j}(w)=
\frac{(-1)^{j-1}}{j-1}h^{l}_{i}(\eta(w'a_{j})\prod_{i=j+1}^{n}a_{i})$$
But by lemma \ref{L:etaneta} again, this is $\frac{(-1)^{j-1}\cdot
j}{j-1}h^{l}_{i}(w)$. Subtraction gives equality.
\end{proof}

\section{Proof of the Main Theorem}

\subsection{Outline of Proof}

Our basic setup is as follows.

\begin{equation} \label{E:CDdpc}
\begin{CD}
\Dpc   @>\sigma_{s_{i}}\sigma_{s_{j}}>> ASS(M(\p)) @>\rho>> Sw(\Dpc) @>\psi_{t}>> Sw(\Dpc)/H' \\
@VV\{\text{$AS,IHX,STU$}\}V  &&  @AfAA @AgAA \\
\Apc && && ASS(\p) @>\psi_{w}>> \W'(p) \\
\end{CD}
\end{equation}

Our maps are as follows:
\begin{enumerate}
    \item  $\sigma_{s_{i}}\sigma_{s_{j}}$ maps acyclic graphs to vertebrates as in Section \ref{SS:sw}
    by selecting $s_{1}$ as the head and $s_{2}$ as the tail. This
    is defined since vertices of elements of $\Dpc$ are ordered.
    \item $\rho$ is defined by figure \ref{Fi:Breaktree}.
    \item $\psi_{t}$ is the quotient map by the relations $H'$ on
    trees.
    \item $\psi_{w}$ is the quotient map by the corresponding relations $H'$ on
    words.
    \item $f$ and $g$ are natural embeddings.
\end{enumerate}

The corresponding commutative diagram in the rooted world is

\begin{equation} \label{E:CDdpcl}
\begin{CD}
\Bbbk(M(\p))   @>\sigma_{s_{i}}^{l}>> rt(ASS(M(\p))) @>\rho^{l}>> Sw(\Dpcl) @>\psi_{t}^{l}>> Sw(\Dpcl)/\Hl \\
@VV\{\text{$AS,IHX,STU$}\}V  &&  @Af^{l}AA @Ag^{l}AA \\
\mathcal{L}(\p) && && ASS(\p) @>\psi_{w}^{l}>> \Wl(p) \\
\end{CD}
\end{equation}

The proof of our main theorem proceeds as follows.
\begin{enumerate}
    \item We prove in Section \ref{SS:rho} that the mapping $\rho$ is well-defined.
    \item We prove in Section \ref{SS:headtail} that $\psi_{t}\rho\sigma_{(i,j)}$ is independent of the mapping $\sigma_{(i,j)}$.
    \item We prove in Section \ref{SS:AS} that the kernel of the mapping $\psi_{w}g$ includes the
    kernels of the $AS$ and of the $IHX$ actions (thus $Sw(\Dpc)/H'\supseteq\Apc$).
\end{enumerate}

As the $H'$ relations come from $AS$ and $IHX$ (thus
$Sw(\Apc)/H'\subseteq\Apc$, this is sufficient to prove isomorphism
between $\Apc$ and $\W'$. The proof in the rooted world is fully
analogous.

Corollaries to the proof, including a new proof to Theorem
\ref{T:hn} are given in Section \ref{SS:hn}.

\subsection{Trees to Swings--- That the $\rho$ Map is Well-Defined}
\label{SS:rho}

The process of breaking down acyclic Jacobi diagrams into sums of
swings over $\Bbbk$ defines a mapping $\rho$ from the space of twice
pointed acyclic Jacobi diagrams $ASS(M(\p))$ to the space of swings,
$Sw(\Dpc)$. The aim of this section is to show that the mapping
$\rho$ (and its analogous mapping $\rho^{l}$ from $\mathcal{L}(\p)$
to $Sw(\mathcal{L}(\p))$) is well-defined--- that it does not depend
on the order in which we break the tree into swings.\par

\begin{lem} \label{L:specialIHX}
The $\rho$ mapping is well defined.
\end{lem}

\begin{proof}
In order to prove this statement, we have to first show that
$\rho$ is independent of the order in which we break down the
branches of a tree until we get a swing.\par

We begin with a twice pointed acyclic Jacobi diagram. For every
trivalent vertex not on the vertebral column and which is not the
neighbour of two legs, we assign two things--- a number to say when
it is to be broken down; and a choice of arc adjacent to the vertex
connecting it to another trivalent vertex `further away' from the
vertebral column. This indicates which subtree is to be broken into
which other subtree. An example is given in Figure 3.

\begin{figure} \label{Fi:rhorder}
\begin{center}
\psfrag{1}{$1$}\psfrag{2}{$2$}\psfrag{3}{$3$}\psfrag{4}{$4$}\psfrag{5}{$5$}
\scalebox{.95}{\includegraphics{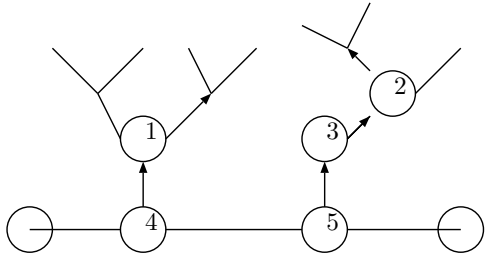}} \caption{A possible
ordering for the breakdown of a tree for the $\rho$ mapping}
\end{center}
\end{figure}

The claim that $\rho$ is independent of this labeling is the claim
that any two such labelings give the same breakdown of the tree into
a sum of swings over $\Bbbk$. It is sufficient to prove this for
every adjacent pair of trivalent vertices.

\begin{figure} \label{Fi:hcVoidIHX}
\begin{center}
\psfrag{t1}{$t_{1}$}\psfrag{t2}{$t_{2}$}\psfrag{t3}{$t_{3}$}\psfrag{t4}{$t_{4}$}
\scalebox{.72}{\includegraphics{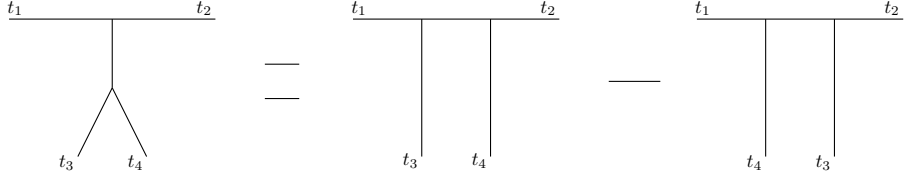}}
\caption{Independence of labeling on two adjacent vertices.}
\end{center}
\end{figure}

Diagrammatically, the claim we have to prove is the claim of Figure
4. Allowing $a_{i}$'s to signify subtrees as well as individual
legs, this is equivalent to Lemma \ref{L:Y2I}.
\end{proof}

This immediately implies that $\rho^{l}$ is also well-defined.

\subsection{Heads and Tails--- Independence from $\sigma_{i,j}$
Map}\label{SS:headtail}

We now prove that it does not matter which legs we chose to be the
head and the tail of our tree before we break down our tree into a
sum of swings--- all such choices are equivalent modulo the $\Hl$
relations.

\begin{lem} \label{L:headindependence}
For any $w_{1},w_{2}$ words in $ASS(\p)$
$$w_{1}\eta(w_{2})=(-1)^{n-1}w_{2}\eta(w_{1})$$
as elements of $\Wl$, where $n:=len(w_{1})+len(w_{2})$.
\end{lem}

\begin{proof}
Let $n_{1}$ be the length of $w_{1}$, $n_{2}$ the length of
$w_{2}$. In order to make the following calculations easier to
understand, let us use the notation
$w_{1}:=\prod_{i=1}^{n_{1}}a_{i}$ and
$w_{2}:=\prod_{i=1}^{n_{2}}b_{i}$.\par

We start with $w_{1}\eta(w_{2})$, applying the $h^{l}_{n}$ action
to all words ending in $b_{n_{2}}$ (which must be of the form
$-w_{1}\eta(\prod_{i=1}^{n_{2}-1}b_{n_{2}-i})b_{n_{2}}$)and the
$h^{l}_{n_{1}+1}$ action to the remaining words (which must be of
the form $w_{1}b_{n_{2}}\eta(\prod_{i=1}^{n_{2}-1}b_{n_{2}-i})$)
and see what happens.

\begin{enumerate}
\item[\textbf{Step 1}] The elements of the result of the action we
have chosen, for which $\eta(w_{2})$ is not split, are
$-(-1)^{n-1}b_{n_{2}}(-1)^{n_{2}}\eta(\prod_{i=1}^{n_{2}-1}b_{i})\eta(w_{1})$,
which by the definition of $h^{l}_{n_{2}}$ is simply
$(-1)^{n-1}w_{2}\eta(w_{1})$.

\item[\textbf{Step 2}]It remains to show that what is left of
$w_{1}\eta(w_{2})$ is zero. Let now $w'_{2}w''_{2}$ be a word in
$\eta(w_{2})$, with $len(w'_{2})>1$
$w'_{2}:=\prod_{i=1}^{len(w'_{2})}b'_{i}$,
$w''_{2}:=\prod_{i=1}^{len(w''_{2})}b''_{i}$. The word
$w'_{2}\eta(w_{1})\tau(w''_{2})$ in the image of our action comes
from exactly two sources. The first is
$h^{l}_{n}(w_{1}w''_{2}\tau(w'_{2}))$ which contains the linear
combination $(-1)^{n-1+len(w''_{2})}w'_{2}\eta(w_{1})\tau(w''_{2})$.
The second place it comes from is $h^{l}_{n}$ of the linear
combination
$(-1)^{len(w''_{2})-1}w_{1}b'_{len(w'_{2})}\tau(w''_{2})\prod_{i=1}^{len(w'_{2}-1)}b'_{i}$
in $w_{1}\eta(\tau(w_{2}))$. The two terms are of opposite sign, and
therefore they cancel each other out.

\item[\textbf{Step 3}] Let us now see what happens when
$len(w'_{2})=1$. Let $w'_{2}:=b'$. $w_{1}\eta(w''_{2})b'$ all give
$(-1)^{n-1+len(w''_{2})}b'\eta(w_{1})\eta(w''_{2})=(-1)^{n_{1}}b'\eta(w_{1})\eta(w''_{2})$
under the $h^{l}_{n}$ action. Taking now
$-\eta(w_{1})b'\eta(w''_{2})$ and acting on it with
$h^{l}_{n_{1}+1}$, we get $(-1)^{n_{1}-1}$ times the same thing.
This is an exact elimination, and it kills all
$b'\eta(w_{1})w''_{2}$ words.\par
\end{enumerate}

This exhausts the ``remainder'', and our lemma is proved.
\end{proof}

For the proof of our main theorem, this lemma is a strong enough
result--- however we may strengthen it still further.

\begin{lem} \label{L:genheadindep}
For any $w_{1},w_{2},\ldots,w_{m}$ words in $ASS(\p)$
$$w_{1}\prod_{i=2}^{m}\eta(w_{i})=(-1)^{m+n-1}w_{m}\eta(\prod_{i=1}^{m-1}(\eta(w_{i})))$$
as elements of $\Wl$, where $n:=\sum_{i=1}^{m}len(w_{i})$. Here
the notation $(\eta(w_{i}))$ means that a function applied to the
word `reads' $\eta(w_{i})$ as if it were a single letter.
\end{lem}

\begin{proof}
We prove the claim by induction. For $m=2$, this is exactly Lemma
\ref{L:headindependence}. Let us assume that the claim holds until
$m=M\in\mathbb{N}$. We proceed as in the proof of that lemma, with
$w_{M+1}$ here playing the the role of $w_{2}$, and the proof is
exactly analogous.
\end{proof}

We now translate Lemma \ref{L:headindependence} to acyclic Jacobi
diagrams.

\begin{cor} \label{C:halfswingindep}
When taking a rooted tree to its left-normed bracketed form by
breaking down trees, the rooted tree we get is independent of which
leg we choose to be the tail of the swing during the construction,
\textbf{as an element of $\Wl$}.
\end{cor}

\begin{proof}

This follows directly from Lemma ~\ref{L:headindependence} (see
Figure 5). On the left hand side, breaking $t_{2}$ into $t_{1}$
gives $(t_{1}\eta(t_{2}))$, while breaking $t_{1}$ into $t_{2}$
gives $(t_{2}\eta(t_{1}))$. By the lemma, these are then equal.

\begin{figure} \label{Fi:headindependence}
\begin{center}
\psfrag{t1}{$t_{1}$}\psfrag{t2}{$t_{2}$}\psfrag{t3}{$t_{3}$}
    \scalebox{.77}{\includegraphics{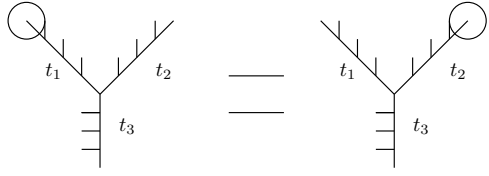}}
    \caption{A graphical interpretation of Lemma ~\ref{L:headindependence}}
\end{center}
\end{figure}

\end{proof}

\begin{cor} \label{C:headindependence}
When ``breaking down a tree into swings'' as in Figure 1, the sum of
swings we get is independent of which legs we choose to be the heads
of the swings during the construction, \textbf{as an element of
$\W'$}.
\end{cor}

\begin{proof}
Let $(b_{1},b_{2})$ and $(b'_{1},b'_{2})$ be two choices of tail and
head respectively for a given tree $T$. Because a tree is
1-connected these exists an arc $c$ which if removed separates the
tree into two subtrees, each of which contains a pair of these four
points. If there exists such an arc $c$ separating the tree into
subtrees containing $(b_{1},b'_{1})$ and $(b_{2},b'_{2})$ in
distinct connected components of $T\setminus\{c\}$, then by
Corollary \ref{C:halfswingindep} the breakdown of these subtrees
into sums of half-swings with heads $b_{1}$ and $b'_{1}$ ($b_{2}$
and $b'_{2}$ respectively) gives the same pre-image under $f$ as
elements of $\Wl$. `Remembering' $c$ proves the corollary for this
case, as the $H$--actions in each $\Wl$ are also in particular
actions of $\W'(deg(T))$. Because of the action $\tau$, the case in
which there exists an arc $c$ separating the tree into subtrees
containing $(b_{1},b'_{2})$ and $(b_{2},b'_{1})$ in distinct
connected components of $T\setminus\{c\}$ is analogous.\par

In the case that there are no such arcs, let $c_{1}$ and $c_{2}$
separate $b_{1}$ and $b_{2}$, $b'_{1}$ and $b'_{2}$ into separate
components respectively. Let us now pick new leaves $b''_{1}$ in
the connected component of $b_{1}$ and $b''_{2}$ in the connected
component $b'_{2}$ of $T\setminus\{c_{1},c_{2}\}$. By
`remembering' $c_{1}$ and $c_{2}$ in turn, the choice of heads
$(b''_{1},b''_{2})$ is equal as an element $\W'$ both to the
choice $(b_{1},b_{2})$ (as now $c_{1}$ separates the tree into
subtrees containing $(b_{1},b''_{1})$ and $(b_{2},b''_{2})$
respectively) and to the choice $(b'_{1},b'_{2})$ (same as before
except with $c_{2}$). Therefore these choices give the same sum of
swings as an element of $\W'$.
\end{proof}

\begin{figure} \label{Fi:Choosehead}
\begin{center}
\psfrag{A}{$A$}\psfrag{B}{$B$}\psfrag{T1}{$T_{1}$}\psfrag{T2}{$T_{2}$}\psfrag{Tm}{$T_{m}$}\psfrag{are
trees.}{are trees.}
    \scalebox{.73}{\includegraphics{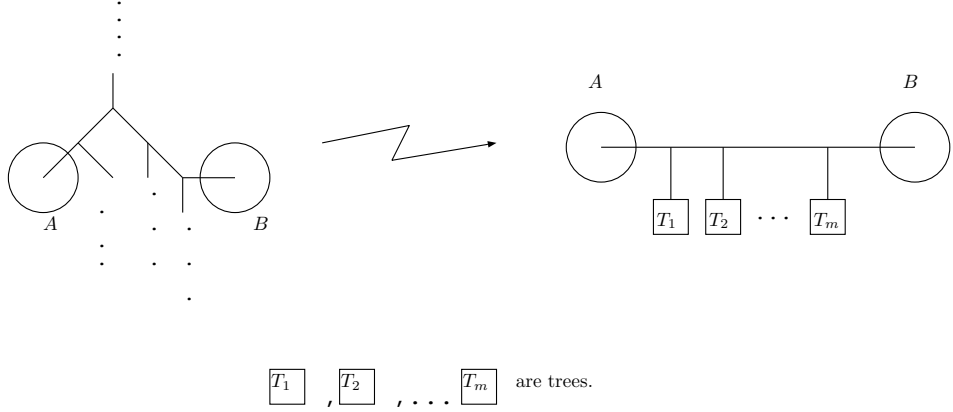}}
    \caption{Choosing a head and a tail, and laying out the graph according to them}
\end{center}
\end{figure}

We can now show the following.

\begin{prop} \label{P:VoidIHX}
The kernel of the $IHX$ relation on $\Apc$, under the $\rho$
mapping, is included in the kernel of the $H$ relations under the
$g$ mapping.
\end{prop}

\begin{proof}
Again, it is enough to prove this for a single tree. We must show
that $IHX$'s which are not in $H'$ do not `impose extra relations
on words'. Let us take such a relation, $IHX_{a}t=t'$.

First, by Lemma \ref{L:headindependence}, the relation is not
dependent on choices of head and tail. Let us take tails of $t_{1}$
and of $t_{2}$ in their left-normed bracketed forms to be our tail
and head respectively. Our relation then takes the form of Figure
4 in Section \ref{SS:rho}, where the proof of Lemma
\ref{L:specialIHX} shows us that this is no new relation.
\end{proof}

We deduce the following.

\begin{cor}
When ``breaking down a tree into swings'' the element in $ASS(\p)$
corresponding with the sum of swings we get is independent of how we
decide to break the tree down, \textbf{as an element of $\W'$}. The
corresponding statement also holds for half-swings.
\end{cor}

\begin{proof}
This is just a combination of Lemma \ref{L:headindependence},
Corollary \ref{C:halfswingindep}, and the Lemma \ref{L:Y2I} as it
is used in the proof of Proposition \ref{P:VoidIHX}.
\end{proof}

\subsection{The $AS$ Relations}\label{SS:AS}

In order to prove the main theorem, it remains only to show that
the kernel of the $AS$ action on $\Dpc$ (and on $M(\p)$) is
contained in the kernel of the $H'$ relations (the $\Hl$
relations).\par

\begin{figure} \label{Fi:Yrel}
\begin{center}
\psfrag{AS}{$AS$}\psfrag{t1}{$t_{1}$}\psfrag{t2}{$t_{2}$}\psfrag{0}{$0$}
\includegraphics{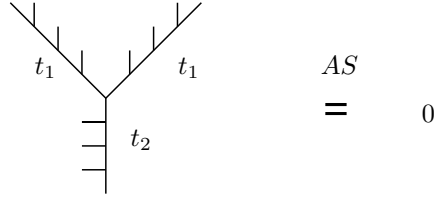}
\caption{The relation $\Y_{t_{1}}$}
\end{center}
\end{figure}

First, we move the problem to the level of words, by defining what
we shall call the $\Y$ relations. We would like the manifestation of
these relations as Jacobi diagram relationships to be as pictured in
Figure 7.

\begin{defn} \label{D:ywvdef}
For $v\in ASS(\p)$, $w$ a word in $ASS(\p)$, let
\begin{equation}
y_{w}(v):=
\begin{cases}
0  &\text{if $v=2w\eta(w)$;}\\
v &\text{otherwise.}
\end{cases}
\end{equation}
\end{defn}

\begin{rem}
The coefficient $2$ in the definition of $y_{w}(v)$ is there to
remind us that the $AS$ relationships have no kernel over a field
of characteristic 2. This coefficient will be ignored from now on,
as long as this point is kept in mind.
\end{rem}

We would like to find a minimal family of such relations which act
on $\Wl$. First, we may demand that $w\in\Wl$. Secondly we may
demand that $len(w)$ be even, as if it is odd and $y_{w}(w')=0$
then $w'=0$ by Lemma \ref{L:headindependence} anyway.

\begin{defn}
Let $\Y$ be the union of all $y_{w}$'s, $w\in\Wl$ a word whose
length is an odd number.
\end{defn}

In order to show that the $\Y$ relations are in the kernel of
$H'$, we shall need the following little lemma.

\begin{lem}[A version of \cite{Reu93}, Theorem 1.4 (v)] \label{L:etanw}
For all $w\in\Wl$ of length $n$, $$\eta(w)=(-1)^{n-1}nw$$.
\end{lem}

\begin{proof}
For $n=2$ we have $\eta(12)=21-12=-2\cdot12$ Let
$w:=\prod_{i=1}^{n}a_{i}$. By induction, using Lemma
\ref{L:etaneta}, we get
$$\eta(wa_{n+1})=(-1)^{n-1}n\cdot wa_{n+1}+a_{n+1}\eta(w)$$ The
second term, by the definition of $h^{l}_{n}$, is equal to
$(-1)^{n-1}wa_{n+1}$, and induction finishes.
\end{proof}

Before we move on, let us point out two pretty little corollaries
to this lemma.

\begin{cor}
If $w_{1}, w_{2}\in\Wl$, and if $\Bbbk$ is not of characteristic
$len(w_{1})$, then $\eta(w_{1})=\eta(w_{2})$ if and only if
$w_{1}=w_{2}$.
\end{cor}

\begin{cor} \label{C:maxlen}
For $\Bbbk$ a field of finite characteristic $n$, there are no
words in $\Wl$ of length greater than $n+1$.
\end{cor}

\begin{proof}
Use of $h^{l}_{n+2}$ twice.
\end{proof}

This turns out to be enough to guarantee that the $\Y$ relations
give us nothing new.

\begin{prop} \label{P:ysize}
$\Y$ is a set of trivial relations on $\Wl$ (and therefore for
$\W'$).
\end{prop}

\begin{proof}
Let $w$ be a word of length $n'$, and let $y_{w}$ be applied to a
word of length $n$. If $\Bbbk$ is of characteristic $n'$, then by
Corollary \ref{C:maxlen}, the claim is trivial. For $n<2n'$ also,
the claim is trivial. Otherwise, for $n\geq2n'$, we have by Lemma
\ref{L:etanw}
$$w\eta(w)w'=\frac{(-1)^{n'-1}}{n'}\eta(w)\eta(w)w'$$

This in turn equals
$\frac{(-1)^{n'}}{n'(2n')}\eta(\eta(w)\eta(w))w'$. By Corollary
\ref{C:etaetakill}, this is zero.
\end{proof}

With this proposition, we have at last completed the proof of our
main theorem, and as promised we have:

\begin{mthm}
$\W'(p)$ is isomorphic to the space $\Apc$, and $\Wl(p)$ is
isomorphic to the space $\mathcal{L}(\p)$.
\end{mthm}

\subsection{Corollaries of the Proof} \label{SS:hn}

As the following proposition shows, this means that there is no
point in picking a head, ``playing around'' and returning to the
same head again in the hope of getting a non-trivial relationship.
This looks fairly obvious, but we could find no proof for it in the
literature.

\begin{lem} \label{L:dblheadchoice}
Let $w:=\prod_{i=1}^{n}a_{i}$ be a word in $ASS(\p)$. Then
$h^{l}_{m}(w)$ ($m\leq n$) followed by choosing the $i$th letter
($i<n$) to be the head of each of the words in $h^{l}_{m}(w)$ by
means of the appropriate fold moves is the same as $h^{l}_{i}(w)$.
\end{lem}

\begin{proof}
For $i=n-1$, this is the same as the proof of lemma
\ref{L:headindependence}.

Now let $c_{a_{i}}(w)$ be the action on an element of $ASS(\p)$ of
choosing $a_{i}$ to be the head. Then
$$h^{l}_{i}(w)=c_{a_{i}}(h^{l}_{i+1}(w))=\cdots=c_{a_{i}}c_{a_{i+1}}\cdots
c_{a_{n-1}}h^{l}_{n}(w)$$ Collapsing this again gives
$c_{a_{i}}h^{l}_{n}(w)$.
\end{proof}

\begin{cor} \label{C:kfunique}
The presentation of an element of $ASS(\p)$ with a given letter as
the head of each word is unique under the action of $\Hl$ (and
therefore under the action of $H'$), \textit{i.e.} if we take a word
$w\in ASS(\underline{p})$ and act on it by arbitrary $h^{l}$--moves,
and then choose $a_{1}$ to be the head of each summand by
$h^{l}$--moves, we recover $w$.

\end{cor}

\begin{proof}
For a single word, this is just Lemma \ref{L:hreduce} and Lemma
\ref{L:dblheadchoice} taken together.
\end{proof}

A further corollary is a proof that the space of acyclic Jacobi
diagrams is isomorphic to the kernel of a mapping between free Lie
algebras.

\begin{proof}[Proof of Theorem \ref{T:hn}]

By the main theorem we identify:

\begin{align*}
\Apc^{n} &\simeq ASS_{n}(\p)/H'\\
\mathcal{L}_{n-1}(\p)\otimes\p &\simeq ASS_{n}(\p)/\Hl-\set{h^{l}_{n}}\\
\mathcal{L}_{n}(\p)&\simeq ASS_{n}(\p)/\Hl
\end{align*}

\noindent where $ASS_{n}(\p)$ denotes the span over $\Bbbk$ of words
of length $n$ over $p$ generators.\par

We begin by defining a map
$g\co\Apc^{n}\To\mathcal{L}_{n-1}(\p)\otimes\p$ by the equation
$g(w)\ass g'(w)-g'(h^{l}_{n}(w))$ where the map $g'$ is defined by
$g'(b_{1}\cdots b_{n})\mapsto b_{2}\cdots b_{n}\otimes b_{1}$ where
the leg labeled $b_{1}$ becomes the root (this is easily seen to be
well--defined by assuming for instance that $w$ is lexically minimal
in its equivalence class)\footnote{This map coincides with the one
used in \cite{HP03,Lev02,Lev06}.}. The mapping
$\mathcal{L}_{n-1}(\p)\otimes\p\To \mathcal{L}_{n}(\p)$ is given by
$a\otimes b\mapsto ba\in\mathcal{L}_{n}(\p)$, a map which we denote
$\ell$.
\par



We first note that $\Apc\subseteq\mathbf{h}_{n}(\p)$--- in other
words that $g(w)=0$ as an element of $\Wl$ for all $w\in\Apc$--- by
choosing $a_{1}$ as the head of each word and applying Corollary
\ref{C:kfunique}.\par

Next, we show the opposite inclusion $\mathbf{h}_{n}(\p)\subseteq
\Apc$.
By the main
theorem, $\ker\ell$ is generated by all $h_{n}^{l}(w)-w$ for all
$n\geq 2$ and for all $w\in ASS(\p)$. To prove that $\mathrm{Im} g$
generates $\ker\ell$, we must show that $\ker\ell$ is generated by
all $h_{n}^{l}(w)-w$ for all $n\geq 2$ and for $w$ belonging to a
maximal set of representatives of equivalence classes in $\Apc^{n}
\simeq ASS_{n}(\p)/H'$. In other words it is sufficient to show that
$h_{n}^{l}(h'_{m}(w))-h'_{n}(w)=h_{n}^{l}(w)-w$ for all $w\in
\mathcal{L}_{n-1}(\p)\otimes\p$.\par

For $m<n$ this follows from Lemma \ref{L:hreduce}. For $m=n$:

\begin{multline*}
h_{n}^{l}(h'_{m}(w))-h'_{n}(w)=(-1)^{2n-1}(a_{1}a_{2}\eta(a_{n}\cdots
a_{3})-a_{1}\eta(a_{n}\cdots a_{3})a_{2})-h'_{n}(w)=\cdots=\\
=\sum_{i=2}^{n-1}(-1)^{n+i-1}a_{1}\cdots a_{i-1}\eta(a_{n}\cdots
a_{i+1})a_{i}\ + (-1)^{2n-1}w-h'_{n}(w)\overset{\text{Lemma
\ref{L:headindependence}}}{=}\\=\sum_{i=2}^{n-1}(-1)^{i-1}a_{n}\cdots
a_{i+1}\eta(a_{1}\cdots a_{i-1})a_{i}\ -h'_{n}(w)-w
=h_{n}^{l}(w)-w
\end{multline*}


Finally as in \cite{Lev02}, note that if we define
$\tilde{g}\co\mathcal{L}(\p)\otimes\p\To\Apc$ by $a\otimes b\mapsto
ba\in\Apc$, we find that

\begin{multline*}\tilde{g}g(w)=\tilde{g}(a_{2}\cdots a_{n}\otimes a_{1}-
h^{l}_{n}(a_{2}\cdots a_{n}\otimes a_{1}))=
w+(-1)^{n}a_{n}\eta(a_{1}\cdots
a_{n-1})=\\w+(-1)^{n-1}a_{n}\eta(a_{1}\cdots
a_{n-2})a_{n-1}+(-1)^{n}a_{n}a_{n-1}\eta(a_{1}\cdots a_{n-2})=\\
2\cdot w+(-1)^{n}a_{n}a_{n-1}a_{n-2}\eta(a_{1}\cdots
a_{n-3})+(-1)^{n-1}a_{n}a_{n-1}\eta(a_{1}\cdots
a_{n-3})a_{n-2}=\cdots=nw\end{multline*}

\noindent where the fourth equality is by $h^{l}_{n-1}$.
\end{proof}

\section{Sample Calculation}\label{S:sample}

In the present section we illustrate the use of our main theorem to
facilitate the calculation by hand of a basis of swings for
$\mathbf{h}_{9}(\underline{9})$. By Witt's dimension formula, the
dimension of this space is:

\begin{multline*}
\dim(\mathbf{h}_{9}(\underline{9}))=\dim(\mathcal{L}_{8}(\underline{9})\otimes
\underline{9})-\dim(\mathcal{L}_{9}(\underline{9}))=\frac{9}{8}(9^{8}-9^{4})-\frac{1}{9}(9^{9}-9^{3})=\\
48420180-43046640=5373540
\end{multline*}

We list the basis elements by their multi-degree. To simplify
notation, in the present section we identify multi-degrees of the
form $(n_{\sigma(1)},n_{\sigma(2)},\ldots,n_{\sigma(9)})$ for all
permutations $\sigma\in \Sigma_{9}$ of the set with nine elements
$\underline{9}$. We also suppress zeros in our notation.\par

We begin by considering the cases where two of the letters appear
only once. In this case, we choose the first pair of such letters as
the head and the tail of the swing, and count the elements which
result.
\begin{itemize}
    \item $\dim(\mathbf{h}(\underline{9})_{(1,1,1,1,1,1,1,1,1)})=7!=5040$
    \item $\dim(\mathbf{h}(\underline{9})_{(1,1,1,1,1,1,1,2)})=9\times8\times\binom{7}{2}\times 5!=181440$
    \item $\dim(\mathbf{h}(\underline{9})_{(1,1,1,1,1,1,3)}=\binom{9}{7}\times 7\times\binom{7}{3}\times 4!=211680$
    \item $\dim(\mathbf{h}(\underline{9})_{(1,1,1,1,1,4)}=\binom{9}{6}\times6\times\binom{7}{4}\times3!=105840$
    \item $\dim(\mathbf{h}(\underline{9})_{(1,1,1,1,5)}=\binom{9}{5}\times5\times\binom{7}{5}\times2!=26460$
    \item $\dim(\mathbf{h}(\underline{9})_{(1,1,1,6)}=\binom{9}{4}\times4\times7\times1=3528$
    \item $\dim(\mathbf{h}(\underline{9})_{(1,1,7)}=\binom{9}{3}\times 3=252$
    \item $\dim(\mathbf{h}(\underline{9})_{(1,1,1,1,1,2,2)}=\binom{9}{7}\times\binom{7}{2}\times\binom{7}{2}\times\binom{5}{2}\times3!=952560$
    \item $\dim(\mathbf{h}(\underline{9})_{(1,1,1,1,2,3)}=\binom{9}{6}\times6\times5\times\binom{7}{3}\times\binom{4}{2}\times2!=1058400$
    \item $\dim(\mathbf{h}(\underline{9})_{(1,1,1,2,4)}=\binom{9}{5}\times5\times 4\times\binom{7}{4}\times 3=264600$
    \item $\dim(\mathbf{h}(\underline{9})_{(1,1,2,5)}=\binom{9}{4}\times4\times 3\times\binom{7}{2}=31752$
    \item $\dim(\mathbf{h}(\underline{9})_{(1,1,1,2,2,2)}=\binom{9}{6}\times\binom{6}{3}\times\binom{7}{2}\times\binom{5}{2}=1058400$
    \item $\dim(\mathbf{h}(\underline{9})_{(1,1,2,2,3)}=\binom{9}{5}\times5\times\binom{4}{2}\times\binom{7}{2}\times\binom{5}{2}=793800$
    \item $\dim(\mathbf{h}(\underline{9})_{(1,1,1,3,3)}=\binom{9}{5}\times\binom{5}{2}\times\binom{7}{3}\times4=176400$
    \item $\dim(\mathbf{h}(\underline{9})_{(1,1,3,4)}=\binom{9}{4}\times 12\times\binom{7}{3}=52920$
\end{itemize}

Of the $450468$ basis elements remaining we next consider those
where only one letter appears only once, which we choose to be the
head. We are then left with subspaces of free Lie algebras which are
homogenous with respect to multi-degree.

\begin{itemize}
\item $\dim(\mathbf{h}(\underline{9})_{(1,2,2,2,2)})=\binom{9}{5}\times
5\times312=196560$
\item $\dim(\mathbf{h}(\underline{9})_{(1,2,2,4)})=\binom{9}{4}\times4\times3\times51=77112$
\item $\dim(\mathbf{h}(\underline{9})_{(1,2,6)})=\binom{9}{3}\times3!\times3=1512$
\item $\dim(\mathbf{h}(\underline{9})_{(1,2,3,3)})=\binom{9}{4}\times4\times 3\times70=105840$
\item $\dim(\mathbf{h}(\underline{9})_{(1,3,5)})=\binom{9}{3}\times3!\times7=3528$
\item $\dim(\mathbf{h}(\underline{9})_{(1,4,4)})=\binom{9}{3}\times3\times8=2016$
\end{itemize}

The basis elements for the last two homogenous subspaces of free Lie
algebras appearing above are those words generated by two letters
$1$ and $2$ in which there is no consecutive sequence an even number
of $2$'s (this is easily calculated by hand using the $\Hl$ moves).
We ask the following question:

\begin{question}
Do words in letters $1$ and $2$ in which there is no consecutive
sequence of an even number of $2$'s form a basis for the free Lie
algebra over two generators?
\end{question}

This question is implicit in \cite{Mos03}, where it is proved for a
few special cases. In the same way we attain the basis elements of
the other homogenous subspaces of free Lie algebras above (`no
sequence of an even number of larger letters between two smaller
letters') suggesting that there may perhaps be something more
general going on.\par

We next look at those $63900$ basis elements in which no letters
appear once, and there is a letter which appears twice. We choose
this letter as head and tail, and then only reflect.

\begin{itemize}
\item $\dim(\mathbf{h}(\underline{9})_{(2,2,2,3)})$ is
$\binom{9}{4}\times4=504$ times
$\frac{1}{2}(\binom{7}{2}\times\binom{5}{2}-3!)=102$ which is
$51408$.
\item $\dim(\mathbf{h}(\underline{9})_{(2,2,5)})$ is
$\binom{9}{3}\times3=252$ times $\frac{1}{2}(\binom{7}{2}-3)=9$
which is $2268$.
\item $\dim(\mathbf{h}(\underline{9})_{(2,3,4)})$ is
$\binom{9}{3}\times 6=504$ times $\frac{1}{2}(\binom{7}{3}-3)=16$
which is $8064$.
\end{itemize}

The final $2160$ cases we deal with individually.

\begin{itemize}
\item $\dim(\mathbf{h}(\underline{9})_{(3,3,3)})$ is
$\binom{9}{3}=84$ times
    \begin{enumerate}
    \item Words beginning with $12$ and ending with $21$. If the
    third letter is one there are four of these. Otherwise there are
    three. Altogether ten possibilities.
    \item Words beginning with $12$ and ending with $31$. If the
    seventh letter is three, then if
    third letter is one we have one possibility, if two there are
    six, and if three then there are three possibilities, altogether
    ten. Otherwise if the seventh letter is one we have six
    possibilities, altogether sixteen.
    \item Words beginning with $13$ and ending with $31$ must be of
    the form $131222331$.
    \end{enumerate}
Altogether 24 possibilities for $84\times24=2016$.
\item $\dim(\mathbf{h}(\underline{9})_{(3,6)})$ is
$\binom{9}{2}\times2=72$ times one, which is $72$.
\item $\dim(\mathbf{h}(\underline{9})_{(4,5)})$ is
$\binom{9}{2}\times2=72$ times a-priori four words of the forms
$121122221$, $121212221$, $121221221$, $122112221$ which can all be
shown to be equivalent modulo the action of $H'$. In total we have
72 possibilities.
\end{itemize}
\bibliographystyle{amsplain}
\bibliography{BigBib}
\end{document}